\documentclass[oneside,american,english,intlimits]{amsart}
\usepackage[T1]{fontenc}
\usepackage[latin9]{inputenc}
\usepackage{amsthm}
\usepackage{amstext}
\usepackage{amssymb}

\makeatletter
\numberwithin{equation}{section}
\numberwithin{figure}{section}
\theoremstyle{plain}
\newtheorem{thm}{Theorem}[section]
  \theoremstyle{plain}
  \newtheorem{assumption}[thm]{Assumption}
  \theoremstyle{definition}
  \newtheorem{defn}[thm]{Definition}
 \theoremstyle{definition}
  \newtheorem{example}[thm]{Example}
  \theoremstyle{remark}
  \newtheorem{notation}[thm]{Notation}
  \theoremstyle{remark}
  \newtheorem{rem}[thm]{Remark}
  \theoremstyle{plain}
  \newtheorem{fact}[thm]{Fact}
  \theoremstyle{plain}
  \newtheorem{cor}[thm]{Corollary}
  \theoremstyle{plain}
  \newtheorem{lem}[thm]{Lemma}
  \theoremstyle{plain}
  \newtheorem{prop}[thm]{Proposition}

\numberwithin{equation}{subsection}
\numberwithin{figure}{section}

\usepackage{graphicx}
\usepackage{xcolor}
\definecolor{darkblue}{rgb}{0,0,0.5} 
\usepackage{transparent}
\usepackage[adobe-utopia]{mathdesign}

\makeatother

\usepackage{babel}

\begin{document}

\title{Curvature-direction measures of self-similar sets}

\author{Tilman Johannes Bohl \and Martina Zähle}

\thanks{Both supported by grant DFG ZA 242/5-1. The first author has previously
worked under the name Tilman Johannes Rothe.}

\email{Tilman.Bohl@uni-jena.de, Martina.Zaehle@uni-jena.de}

\address{Friedrich Schiller University Jena}



\maketitle
\newcommand\printornot[1]{#1}

\global\long\def\op#1{\operatorname{#1}}
\global\long\def\Int{\op{int}}
\global\long\def\Frac{\op{frac}}
\global\long\def\esup{\op{ess\, sup}}
\global\long\def\var{\op{var}}
\global\long\def\diam{\op{diam}}
\global\long\def\orth{\operatorname{orth}}
\global\long\def\id{\operatorname{id}}
\global\long\def\dist{\operatorname{dist}}
\global\long\def\assign{\mathrel{:=}}
\global\long\def\G{G}
\global\long\def\H{\mathcal{H}}
\global\long\def\Hl{\H_{\G}}
\global\long\def\ra{\rightarrow}
\global\long\def\ep{\epsilon}
\global\long\def\r{\rho}
\global\long\def\Rd{\mathbb{R}^{d}}
\global\long\def\R{\mathbb{R}}
\global\long\def\N{\mathbb{N}}
\textit{}\global\long\def\Z{\mathbb{Z}}
\global\long\def\borel{\mathcal{B}}
\global\long\def\B{\borel}
\global\long\def\dirs{R}

\begin{abstract}
We obtain fractal Lipschitz-Killing curvature-direction measures for
a large class of self-similar sets $F$ in $\Rd$. Such measures jointly
describe the distribution of normal vectors and localize curvature
by analogues of the higher order mean curvatures of differentiable
sub-manifolds. They decouple as independent products of the unit Hausdorff
measure on $F$ and a self-similar fibre measure on the sphere, which
can be computed by an integral formula. The corresponding local density
approach uses an ergodic dynamical system formed by extending the
code space shift by a subgroup of the orthogonal group. We then give
a remarkably simple proof for the resulting measure version under
minimal assumptions.
\end{abstract}

\section{Introduction}

The \textquotedbl{}second order\textquotedbl{} anisotropic structure
of self-similar sets $F$ in $\Rd$ is studied by means of approximation
with parallel sets $F(\ep)$ of small distances $\epsilon$. This
leads to fractal curvature-direction measures and their local {}``densities''.
From the isotropic point of view this was first investigated in the
pioneering work by Winter \cite{Win08MR2423952} (deterministic self-similar
sets with polyconvex parallel sets, curvature measures) and then in
\cite{Zaehle11SelfsimRandomFractalsMR2763731}(self-similar random
sets with singular parallel sets, total curvatures), \cite{WinterZaehle10cmofsss}
(deterministic measure version for singular parallel sets) and \cite{RatajZaehle10CurvatureDensitiesSelfSimilarSets}
(dynamical approach to curvatures measures and their local densities).
The special case of the Min\-kow\-ski content was treated earlier,
e.g. in \cite{LapidusPomerance93MR1189091} and \cite{Falconer95MinkowskiMeasurabilityMR1224615}
($d=1$), in \cite{Gatzouras00LacunarityMR1694290} (self-similar
random fractals for any $d$), in \cite{KombrinkKessebohmer10OneDimSelfconformal}
($d=1$ self-conformal sets), and in a very general context recently
in \cite{RatajWinter09MeasuresOfParallelSetsArxiv09053279}.\\
 In the present paper we extend these results to anisotropic quantities
for the fractal sets (cf. the remark at the end of the paper). We
mainly follow the dynamical approach from \cite{RatajZaehle10CurvatureDensitiesSelfSimilarSets}
and give a new and short proof for convergence of the corresponding
measures under weaker assumptions, which considerably simplifies the
former approaches.\\
 The classical geometric background are extensions of Federer's
curvature measures for sets of positive reach (\cite{Federer59MR0110078})
regarding normal directions, the so-called curvature-direction measures
or generalized curvature measures (cf. \cite{Schneider80PolyconvexSteinerFormulaMR566443},
\cite{Zaehle86IntRepMR849863}, and various subsequent papers). In
the context of differential geometry they would correspond to the
integrals of higher order mean curvatures over those points of the
boundary which together with the unit normal belong to a fixed Borel
set in $\Rd\times S^{d-1}$. Since the boundaries under consideration
are non-smooth, we are working with integrals of generalized mean
curvatures over their unit normal bundles in the sense of \cite{RatajZaehle05GeneralNormalMR2131910}.
Anisotropic curvature quantities already prove useful to describe
heterogeneous materials in the classical setting, see \cite{SchroederturkMecke11MinkowskiShapeAnalysis}
and the references therein.

Under some regularity condition on the parallel sets, which is always
fulfilled if $d\leq3$ (according to a result in \cite{Fu85TubularNeighMR816398})
or if the convex hull of $F$ is a polytope (\cite{Pokorny11CriticalValuesSelfSimilar1101.1219}),
these signed measures $C_{k}(F(\ep),\cdot)$ on $\Rd\times S^{d-1}$
of order $k$ are determined for almost all $\ep$ (see Section \ref{sub:Curvature-measures-of}).
Let $G$ be the subgroup generated by the orthogonal components of
the similarities associated with $F$. We first restrict the curvature-direction
measures to covariant sets $A_{F}(x,\ep)\times\phi\dirs$ from a suitable
neighborhood net in $\Rd\times S^{d-1}$, e.g. to $(F(\ep)\cap B(x,a\ep))\times\phi\dirs$,
where $a>1$ is fixed, $x\in F$, $\dirs$ is a Borel set in $S^{d-1}$
and $\phi$ an element of $G$. Under some integrability condition
we obtain for any $\dirs$ and almost all $(x,\phi)\in F\times G$
an average limit of the rescaled versions \[
\ep^{-k}C_{k}\big(F(\ep),A_{F}(x,\ep)\times\phi\dirs\big)\]
 as $\ep\rightarrow0$ (Theorem \ref{curvdens}). For the total value
at $\dirs=S^{d-1}$ this corresponds to the \textit{$k$-th fractal
curvature densities} introduced in \cite{RatajZaehle10CurvatureDensitiesSelfSimilarSets}.
As function in $\dirs$ the limits may be considered as second order
fibre measures which, divided by the total values, describe the \textit{local
fractal direction 'distributions' weighted by curvatures}. Due to
self-similarity and the behavior of the classical curvature-direction
measures these limits do not depend on the base points $x$ and on
$\phi$, i.e., we obtain at almost all points constant fractal curvature
densities, constant fibre measures and therefore also constant local
curvature-direction distributions. In this part of our approach an
extended ergodic dynamical system together with the geometric scaling
and invariance properties of $F$, $A_{F}$ and $C_{k}$ provide the
main tools.\\

Then an appropriate choice of the neighborhood net $A_{F}$, similar
as in \cite{RatajZaehle10CurvatureDensitiesSelfSimilarSets}, and
a new method of proof is used in order to derive the corresponding
measure result: Under some slightly stronger uniform integrability
condition the rescaled curvature-direction measures \[
\ep^{D-k}C_{k}(F(\ep),\cdot)\]
 weakly converge in the average as $\ep\rightarrow0$ to a \textit{fractal
curvature-direction measure} which is the product of the normalized
Hausdorff measure $\H^{D}$ on $F$ and the above fibre measure (Theorem
\ref{thm:c_k converge}).\\
 (The fractal curvature measures of order $k=0,1,\ldots,d$ from
\cite{RatajZaehle10CurvatureDensitiesSelfSimilarSets} arise as the
projection measures on the base point component. Because of self-similarity
they are all constant multiples of the normalized Hausdorff measure
on $F$ and the constants agree with the above fractal curvature densities.)\\
 The new approach makes the proof much shorter, even in our non-isotropic
version. It enlightens the essential measure theoretic background
and shows that the uniform integrability condition is sharp.\\

Finally, a modified Sierpinski gasket is discussed as an example.

\section{Basic Notions}

\subsection{Self-similar sets, code space, and measures}

The notion of self-similar sets is well-known from the literature.
See Hutchinson \cite{Hutchinson81FractalsMR625600} for the first
general approach and the relationships mentioned below without a reference.
We use here the following notation and results.\\
 The basic space is a compact set $J\subset\Rd$ with $J=\overline{\Int J}$.
$S_{1},\ldots,S_{N}$ denotes the generating set of \textit{contracting
similarities} in $\Rd$ with \textit{contraction ratios} $r_{1},\ldots,r_{N}$
and orthogonal group components $\varphi_{1},\dots,\varphi_{N}$.

We assume the \textit{strong open set condition} (briefly (SOSC))
with respect to $\Int J$, i.e., \[
\bigcup_{j=1}^{N}S_{j}(J)\subset J,\qquad S_{j}(\Int J)\cap S_{l}(\Int J)=\emptyset\,,\quad j\neq l,\]
and that there exists a sequence of indices $l_{1},l_{2},\ldots,l_{m}\in\{1,\ldots,N\}$
such that \[
S_{l_{1}}\circ S_{l_{2}}\ldots\circ S_{l_{m}}(J)\cap\Int J\neq\emptyset\,.\]
 The latter (strong) condition is here equivalent to \[
F\cap\Int J\ne\emptyset\]
 where $F$ denotes the associated self-similar fractal set $F$.
(According to a result of Schief (see \cite{Schief94SeparationMR1191872}),
(SOSC) for some $J$ is already implied by the open set condition
on the similarities. A characterization of (SOSC) in algebraic terms
of the $S_{i}$ is given in Bandt and Graf \cite{BandtGraf92MR1100644}.)
The set $F$ may be constructed by means of the \textit{code space}
$W:=\{1,\ldots,N\}^{\N}$, the set of all infinite words over the
alphabet $\{1,\ldots,N\}$. We write $W^{n}:=\{1,\ldots,N\}^{n}$
for the set of all words of length $|w|=n$, $W^{*}:=\bigcup_{n=1}^{\infty}W^{n}$
for the set of all finite words, $w|n:=w_{1}w_{2}\ldots w_{n}$ if
$w=w_{1}w_{2}\ldots w_{n},w_{n+1}\ldots$ for the \textit{restriction}
of a (finite or infinite) word to the first $n$ components, and $vw$
for the \textit{concatenation} of a finite word $v$ and a word $w$.
If $w=w_{1}\ldots w_{n}\in W^{n}$ we also use the abbreviations $S_{w}:=S_{w_{1}}\circ\ldots\circ S_{w_{n}}$,
$\varphi_{w}:=\varphi_{w_{1}}\circ\ldots\circ\varphi_{w_{n}}$, and
$r_{w}:=r_{w_{1}}\ldots r_{w_{n}}$ for the contraction ratio of this
mapping. Finally we denote $K_{w}:=S_{w}(K)$ for any compact set
$K$ and $w\in W^{*}$. (For completeness we also write $K_{\emptyset}:=K$.)
In these terms the set $F$ is determined by \[
F=\bigcap_{n=1}^{\infty}\bigcup_{w\in W^{n}}J_{w}\]
 and characterized by the \textit{self-similarity property} $F=S_{1}(F)\cup\ldots\cup S_{N}(F)$.
Iterated applications yield \[
F=\bigcup_{w\in W^{n}}F_{w},\quad n\in\N.\]
 As in the literature, we will use the abbreviation \[
S(K):=\bigcup_{j=1}^{N}S_{j}(K)\]
 for compact sets $K$, i.e., $F=S^{n}(F)$, $n\in\N$.\\
 Alternatively, the self-similar fractal $F$ is the image of the
code space $W$ under the \textit{projection} $\pi$ given by \[
\pi(w):=\lim_{n\rightarrow\infty}S_{w|n}x_{0}\]
 for an arbitrary starting point $x_{0}$. The mapping $w\mapsto x=\pi(w)$
is bi-unique except for a set of points $x$ of $D$-dimensional Hausdorff
measure $\H^{D}$ zero, and the \textit{Hausdorff dimension} $D$
of $F$ is determined by \begin{equation}
\sum_{j=1}^{N}r_{j}^{D}=1\,.\label{dim}\end{equation}
 Up to exceptional points we \textit{identify} $x\in F$ with its
\textit{coding sequence} and write $x_{1}x_{2}\ldots$ for this infinite
word, i.e. $\pi(x_{1}x_{2}\ldots)=x$, and write \[
x|n:=x_{1}\ldots x_{n}\]
 for the corresponding finite words.\\
 If $\nu$ denotes the infinite product measure on $W$ determined
by the probability measure on the alphabet $\{1,\ldots,N\}$ with
single probabilities $r_{1}^{D},\ldots,r_{N}^{D}$, then the normalized
$D$-dimensional Hausdorff measure with support $F$ equals \begin{equation}
\mu_{F}:=\H^{D}(F)^{-1}\H^{D}(F\cap(\cdot))=\nu\circ\pi^{-1}\,.\label{bernoulli-meas}\end{equation}
 It is also called the natural \textit{self-similar measure} on $F$,
since we have \begin{equation}
\mu_{F}=\sum_{j=1}^{N}r_{j}^{D}\mu_{F}\circ S_{j}^{-1}\,.\label{ssmu}\end{equation}

$\G$ is the compact group generated by all the $\varphi_{j}$, i.e.
the orthogonal group components of the $S_{j}$: \begin{equation}
\G\assign\op{cl}\left\{ \varphi_{j}\,:\, j=1,\dots,N\right\} \subseteq O\left(d\right)\label{eq:G formel}\end{equation}
 Denote its normalized Haar measure by $\H_{\G}$.

From (SOSC) on $J$ one obtains \begin{equation}
\int|\ln d(y,J^{c})|\, d\mu_{F}(y)<\infty\label{intlogdist}\end{equation}
 (\cite[Proposition 3.4]{Graf95TangentialMR1363139}) and $\mu_{F}\left(\partial J\right)=0$.
Furthermore, by the open set condition $F$ is a $D$-\textit{set
(Ahlfors regular)}, i.e., there exist positive constants $c_{F}$
and $C_{F}$ such that \begin{equation}
c_{F}\, r^{D}\le\H^{D}(F\cap B(x,r))\le C_{F}\, r^{D}\,,~~x\in F,~r\le\diam F\,.\label{D-set}\end{equation}

Finally, we write $f\otimes g\left(x,n\right)=f\left(x\right)g\left(n\right)$
if $f:\Rd\ra\R$ and $g:S^{-1}\ra\R$, and $\mu\otimes\nu$ for the
product measures of $\mu$ and $\nu$, and $\mu\left(f\right)$ for
$\int fd\mu$.

\subsection{\label{sub:Curvature-measures-of}Curvature-direction measures of
parallel sets}

\label{sec:classcurv} We will use the following notations for points
$x$ and subsets $E$ of $\Rd$: \[
d(x,E):=\inf_{y\in E}|x-y|~,~~|E|:=\diam E=\sup_{x,y\in E}|x-y|\,.\]
 The background from classical singular curvature theory is summarized
in \cite{Zaehle11SelfsimRandomFractalsMR2763731}. We recall some
of those facts. For certain classes of compact sets $K\subset\Rd$
(including many classical geometric sets) it turns out that for Lebesgue-almost
all distances $r>0$ the parallel set ($r$-tube, offset, Minkowski
sausage) \[
K(r):=\left\{ x:d(x,K)\text{\ensuremath{\le}}r\right\} \]
 possesses the property that the closure of its complement \[
\widetilde{K(r)}\assign\overline{K(r)^{c}}\]
 is a set of positive reach in the sense of Federer \cite{Federer59MR0110078}
with Lipschitz boundary. A sufficient condition is that $r$ is a
regular value of the Euclidean distance function to $K$ (see Fu \cite[Theorem 4.1]{Fu85TubularNeighMR816398}
together with \cite[Proposition 3]{RatajZaehle03NormalCyclesByApproxMR1983898}).
(In $\mathbb{R}^{2}$ and $\mathbb{R}^{3}$ this is fulfilled for
all $K$, see \cite{Fu85TubularNeighMR816398}, and in any $\Rd$
for self-similar sets whose convex hull is a polytope, see \cite{Pokorny11CriticalValuesSelfSimilar1101.1219}.)
In this case both the sets $\widetilde{K(r)}$ and $K(r)$ are \textit{Lipschitz
d-manifolds of bounded curvature} in the sense of \cite{RatajZaehle05GeneralNormalMR2131910},
i.e., their \textit{k-th Lipschitz-Killing curvature-direction measures},
$k=0,1,\ldots,d-1$, are determined in this general context and agree
with the classical versions in the special cases. Their marginal $C_{k}\left(K\left(r\right),\Rd\times\text{·}\right)$
is known as area measure in convex geometry. Moreover, they satisfy
\begin{equation}
C_{k}\left(K\left(r\right),\text{·}\right)=\left(-1\right)^{d-1-k}C_{k}\left(\widetilde{K(r)},\rho\left(\cdot\right)\right),\label{eq:sign}\end{equation}
where $\rho\left(x,n\right)=\left(x,-n\right)$ is the reflection
on $\Rd\times S^{d-1}$. The $C_{k}(K(r),\cdot)$ are signed measures
with finite \textit{variation measures} $C_{k}^{\var}(K(r),\cdot)$.
Their explicit integral representations are reduced to \cite{Zaehle86IntRepMR849863}
(cf. \cite[Theorem 3]{RatajZaehle05GeneralNormalMR2131910} for the
general case).

$C_{d-1}(K(r),\cdot\times S^{d-1})$ agrees with one half of the $(d-1)$-\textit{dimensional
Hausdorff measure} $\H^{d-1}$ on the boundary $\partial K(r)$. For
completeness, we define $C_{d}(K(r),\cdot)$ as \textit{Lebesgue measure
restricted to} $K(r)$, times the unique rotation invariant probability
measure on $S^{d-1}$. The \textit{total measures (curvatures)} of
$K(r)$ are denoted by\[
C_{k}\left(K\left(r\right)\right)\assign C_{k}\left(K\left(r\right),\Rd\times S^{d-1}\right),\, k=0,\dots,d.\]
By an associated Gauss-Bonnet theorem (see \cite[Theorems 2,3]{RatajZaehle03NormalCyclesByApproxMR1983898})
the \textit{total Gauss curvature} $C_{0}(K(r))$ coincides with the
\textit{Euler-Poincaré characteristic} $\chi(K(r))$.\\
In the present paper only the following main properties of the
curvature measures for such parallel sets will be used, $A\subseteq\Rd,\dirs\subseteq S^{d-1}\,\mbox{Borel}$:\\
The curvature measures are \textit{motion invariant}, i.e.,\begin{gather}
C_{k}\left(g\left(K\left(r\right)\right),gA\times\varphi_{g}\dirs\right)=C_{k}\left(K\left(r\right),A\times\dirs\right)\label{eq:covariance-curv}\\
\mbox{for any Euclidean motion }g,\,\mbox{orthogonal component }\varphi_{g}\nonumber \end{gather}
 they are \textit{homogeneous of degree} $k$, i.e., \begin{equation}
C_{k}\left(\lambda K\left(r\right),\left(\lambda A\right)\times\dirs\right)=\lambda^{k}C_{k}\left(K\left(r\right),A\times\dirs\right),\quad\lambda>0,\label{eq:homogenity-curv}\end{equation}
and locally determined, i.e., \begin{equation}
C_{k}\left(K\left(r\right),\left(\cdot\right)\cap O\times S^{d-1}\right)=C_{k}\left(K^{\prime}\left(r^{\prime}\right),\left(\cdot\right)\cap O\times S^{d-1}\right)\label{eq:locality-curv}\end{equation}
for any open set $O\subset\Rd$ such that $K\left(r\right)\cap O=K^{\prime}\left(r^{\prime}\right)\cap O$,
where $K(r)$ and $K^{\prime}\left(r^{\prime}\right)$ are both parallel
sets where the closures of the complements have positive reach.

\section{Curvature-direction measures of self-similar sets}

\subsection{Covariant local neighborhood nets}

Throughout this entire paper, we will assume:
\begin{assumption}
\label{ass:Neigborhood-regularity}(Regularity of parallel sets) For
Lebesgue almost all $\epsilon>0$:
\begin{enumerate}
\item reach $\widetilde{F\left(\epsilon\right)}>0$,
\item $\op{nor}\widetilde{F\left(\epsilon\right)}\cap\rho\op{nor}\widetilde{F\left(\epsilon\right)}=\varnothing$.
\end{enumerate}
\end{assumption}
Here $\op{nor}X$ denotes the unit normal bundle of a set with positive
reach (as subset of $\Rd\times S^{d-1}$, cf. the reference papers)
and $\rho$ the normal reflection $\left(x,n\right)\mapsto\left(x,-n\right)$.
Curvature measures of $\widetilde{F\left(\epsilon\right)}$ exist
under these conditions.

Many classes of sets satisfy this requirement, e. g. in case $F\left(\epsilon\right)$
is polyconvex, or any set $F$ in space dimension $d\leq3$ \cite{Fu85TubularNeighMR816398},
or if $\op{conv}F$ is a polytope \cite{Pokorny11CriticalValuesSelfSimilar1101.1219}.
In higher dimensions, it is more convenient to check that Lebesgue
almost all $\epsilon$ are regular values of the Euclidean distance
function \cite[Theorem 4.1]{Fu85TubularNeighMR816398}, \cite[Proposition 3]{RatajZaehle03NormalCyclesByApproxMR1983898}.

We conjecture the assumption is always true: Is is an open problem
whether almost all $\epsilon$ are always regular distances to any
given deterministically self-similar set with the Open Set Condition.
The only known counterexample is not self-similar \cite{Ferry76BoundaryManifoldMR0413112}.

In order to refine the results from \cite{RatajZaehle10CurvatureDensitiesSelfSimilarSets}
and \cite{WinterZaehle10cmofsss} for the anisotropic case we will
consider the covariant local neighborhood nets $\left\{ A_{F}\left(x,\epsilon\right)\subseteq\Rd\times S^{d-1}\,:\, x\in F,\,0<\epsilon<\epsilon_{0}\right\} $
from \cite{RatajZaehle10CurvatureDensitiesSelfSimilarSets} marked
by direction sets $\dirs\in\borel\left(S^{d-1}\right)$. As a main
step we will show that for a certain constant $b$, the measures \begin{equation}
\Delta_{F,x,\delta}\left(\dirs\right)\assign\frac{1}{\left|\ln\delta\right|}\int_{\delta}^{d\left(x,J^{c}\right)/b}\epsilon^{-k}C_{k}\left(F\left(\epsilon\right),\, A_{F}\left(x,\epsilon\right)\times\dirs\right)\epsilon^{-1}d\epsilon\label{eq:area measure of finite parallel set}\end{equation}
 on $S^{d-1}$ converge weakly as $\delta\ra0$ for $\mu_{F}$-almost
all $x$ (Corr. \ref{pro:area measures converge, invariant}). The
limit measure $D_{C_{k}^{\op{frac}}|F}$ is the directional component
of the fractal curvature-direction measure $C_{k}^{\op{frac}}\left(F,\cdot\right)$
to be derived on $F\times S^{d-1}$. To this aim we first consider
the sets $A_{F}\left(x,\epsilon\right)\times\phi\left(\dirs\right)$
for a fixed $\dirs\in\B\left(S^{d-1}\right)$ and calculate the limit\[
\lim_{\delta\ra0}\Delta_{F,x,\delta}\left(\phi\dirs\right)=D_{C_{k}^{\op{frac}}|F}\left(\phi\dirs\right)\]
 for $\mu_{F}$-almost all $x\in F$ and $\H_{\G}$-almost all $\phi\in\G$
(Theorem \ref{thm:area measure}). The main tools are to translate
the problem into the language of an extended shift dynamical system,
and Birkhoff's ergodic theorem.
\begin{defn}
\label{def:covariant neighborhood net}Let constants $a>1$ and $\epsilon_{0}>0$
be given, and denote \[
b\assign\mbox{max}\left(2a,\epsilon_{0}^{-1}\left|J\right|\right).\]
 A \textit{covariant neighborhood net }in $F$ is a family of measurable
sets \[
\left\{ A_{F}\left(x,\epsilon\right)\subseteq\Rd\times S^{d-1}\,:\, x\in F,\,0<\epsilon<\epsilon_{0}\right\} \]
 satisfying the following three conditions.
\begin{enumerate}
\item $A_{F}\left(x,\epsilon\right)\subseteq\left(F\left(\epsilon\right)\cap B\left(x,a\epsilon\right)\right)$,
\item $A_{F}\left(x,\epsilon\right)=S_{j}A_{F}\left(S_{j}^{-1}\left(x\right),r_{j}^{-1}\epsilon\right)$
if $1\leq j\leq N$, $x\in F_{j}$, and $\epsilon<d\left(x,\left(S_{j}J\right)^{c}\right)/b$.\label{enu:loc-hom-neighbor A_F}
\item The indicator function $1_{A_{F}\left(x,\ep\right)}\left(z\right)$
is a measurable function of $\left(x,\ep,z\right)\in F\times\left(0,\epsilon_{0}\right)\times\Rd$.
\end{enumerate}
\end{defn}
(Note that $r_{j}^{-1}\epsilon<\epsilon_{0}$ in \eqref{enu:loc-hom-neighbor A_F}.) 
\begin{example}
Two possible choices of $A_{F}$ are:\begin{eqnarray}
A_{F}(x,\ep) & \assign & \left(F(\ep)\cap B(x,a\ep)\right),\quad\ep>0,\label{eq:A_F in ball}\\
A_{F}(x,\ep) & \assign & \left\{ x^{\prime}\in F(\ep):|x-x^{\prime}|\le\r_{F}(x^{\prime},\ep)\right\} ,\label{eq:A_F useful}\\
 &  & \quad\,\quad\qquad\;\;\;0<\ep<\ep_{0}:=\H^{D}(F)^{1/D}\nonumber \end{eqnarray}
 where $\r_{F}(z,\ep)$ is determined for $0<\ep<\ep_{0}$ by the
condition \[
\rho_{F}(x^{\prime},\ep)=\min\{\rho:\,\H^{D}(F\cap B(x^{\prime},\rho))=\ep^{D}\}.\]
The choice \eqref{eq:A_F useful} is important because of its close
ties with the measure version of our main results, as discussed in
section \ref{sec:Global-versus-Local}, and is made to match \eqref{eq:C_k <-> local version}.
The requirements of the definition are met, see \cite[Lemma 2.1.2]{RatajZaehle10CurvatureDensitiesSelfSimilarSets}.
\end{example}

\subsection{Main result 1: fibre measures}

Recall $\G$ from \eqref{eq:G formel}.
\begin{notation}
For Borel sets $\dirs\subseteq S^{d-1}$, define\begin{eqnarray}
D_{C_{k}^{\Frac}|F}(\dirs) & \assign & \frac{1}{\sum_{j=1}^{N}r_{j}^{D}|\ln r_{j}|}\,\int_{F\times\G\,}\:\int_{d(x,\,(S_{x_{1}}J)^{c})/b}^{d(x,\, J^{c})/b}\label{curvdens-int}\\
 &  & \ep^{-k}C_{k}\big(F(\ep),A_{F}\left(x,\ep\right)\times\phi\left(\dirs\right)\big)\,\,\ep^{-1}d\ep\,\;\: d\left(\mu_{F}\otimes\H_{\G}\right)(x,\phi)\nonumber \end{eqnarray}
 provided the integral exists. It serves as a limit object, and we
will show it is a (signed) measure. $D_{C_{k}^{\Frac}|F}^{\pm}$ is
defined by substituting the variation measure $C_{k}^{\pm}$ for $C_{k}$
in the formula.
\end{notation}
The case of an infinite limit is only relevant to prove our results
are sharp.
\begin{thm}
\label{curvdens}(Fibre measure) Suppose that the self-similar set
$F$ in $\Rd$ with Hausdorff dimension $D$ satisfies the strong
open set condition w.r.t. $\Int J$. Let the system of all $A_{F}(x,\ep)$
with $x\in F$, $\ep<\ep_{0}$, be a covariant neighborhood net with
constants $a>1$ and $\ep_{0}>0$. Let $b=\max\big(2a,\ep_{0}^{-1}|J|)\big)$.
Let $k\in\{0,1,\ldots,d\}$. If $k\le d-2$, we additionally suppose
the regularity of parallel sets: Assumption \ref{ass:Neigborhood-regularity}.
Recall \[
\Delta_{F,x,\delta}\left(\dirs\right)\assign\frac{1}{|\ln\delta|}\int_{\delta}^{d(x,J^{c})/b}\ep^{-k}C_{k}\big(F(\ep),\, A_{F}(x,\ep)\times\dirs)\big)\,\,\ep^{-1}d\ep.\]

Then for any fixed Borel set $\dirs\subseteq S^{d-1}$ and $\mu_{F}\otimes\H_{\G}$-a.a.
$\left(x,\phi\right)\in F\times\G$ we have the following: \begin{equation}
\lim_{\delta\rightarrow0}\Delta_{F,x,\delta}\left(\phi\left(\dirs\right)\right)=D_{C_{k}^{\Frac}|F}(\dirs)\label{dens}\end{equation}
 provided the double integral in \eqref{curvdens-int} converges absolutely.
For $k\in\{d-1,d\}$ it always converges absolutely.

The assertion remains true (even without assuming integrability) if
$C_{k}$ is replaced everywhere with $C_{k}^{+}$, $C_{k}^{-}$, or
$C_{k}^{\var}$, and accordingly $D_{C_{k}^{\Frac}|F}$.\end{thm}
\begin{rem}
For the choice of $A_{F}$ as in \eqref{eq:A_F useful}, the mapping
$\dirs\mapsto D_{C_{k}^{\Frac}|F}(\dirs)$ is the constant fibre of
the disintegration of the associated $k$-th fractal curvature-direction
measure over the Hausdorff measure, see section \ref{sec:Global-versus-Local}.
\end{rem}
The case $k=d$ treats the Minkowski content and is formally included,
but $D_{C_{d}^{\Frac}|F}$ is always isotropic and does not provide
geometric information beyond its total mass. Even that is redundant,
because whenever $D<d$, \[
D_{C_{d}^{\Frac}|F}=\frac{1}{d-D}\, D_{C_{d-1}^{\Frac}|F},\]
 and one side exists whenever the other one does, see \cite[Theorems 4.5, 4.7]{RatajWinter09MeasuresOfParallelSetsArxiv09053279}
and Theorem \ref{thm:c_k converge}.
\begin{rem}
\label{rem:sufficient for absolute convergence}Sufficient (sharper)
conditions for both the integrability of \eqref{curvdens-int} in
the above theorem and the uniform integrability in Theorem \ref{thm:c_k converge}
are
\begin{itemize}
\item \begin{eqnarray}
\underset{0<\ep<\ep_{0},\, x\in F}{\esup}\ep^{-k}\, C_{k}^{\var}\big(F(\ep),B\left(x,a\ep\right)\times S^{d-1}\big) & < & \infty,\,\label{eq:ck_var bounded}\end{eqnarray}
 i.e., the rescaled curvature of almost every $a\ep$-balls is bounded;
or
\item condition \cite[Theorem 2.1(ii)]{WinterZaehle10cmofsss}, i.e, the
above but restricted to {}``overlap sets''; or 
\item polyconvex parallel sets: $F(\ep)$ is a finite union of convex sets
for any and therefore all $\ep>0$; or
\item $k=d$ or $k=d-1$, i.e., Minkowski or surface content.
\end{itemize}
We welcome future work to search for even more convenient conditions.
One can show Winter's Strong Curvature Bound Condition \cite[Theorem 5.1]{Winter10CurvBoundsArXiv10102032}
also belongs on this list.\end{rem}
\begin{proof}
(of the remark) Using \eqref{intlogdist}, one can see the supremum
condition implies (uniform) integrability by the same proof as \cite[Rem 3.1.3]{RatajZaehle10CurvatureDensitiesSelfSimilarSets}.
The remaining sufficient conditions reduce to the first one: For the
overlaps-only condition, adapt the proof of \cite[Lemma 3.3]{WinterZaehle10cmofsss}.
For polyconvex parallel sets, the modifications to the proof of \cite[Lemma 5.3.2]{Win08MR2423952}
are sketched in \cite[Remark 2.2.3]{RatajZaehle10CurvatureDensitiesSelfSimilarSets}.
For the Minkowski or surface content, see \cite[Remark 3.2.2]{RatajZaehle10CurvatureDensitiesSelfSimilarSets}.
\end{proof}
Under the conditions of the remark, there is a more convenient formula
for the total mass of $D_{C_{k}^{\Frac}|F}$, see \cite[Theorem 2.3.6]{Win08MR2423952}
and \cite[Theorem 2.3.8]{Zaehle11SelfsimRandomFractalsMR2763731}.
We conjecture the formula holds much more generally.

Section \ref{sec:Example} gives some examples that use these conditions.

\subsection{\label{sub:Extended-shift-dynamical}Extended shift dynamical system
and proof of the theorem}

As an essential auxiliary tool for the proof we use the ergodic\textit{
dynamical system} $\left[W\times\G,\nu\otimes\H_{\G},\theta\right]$
on the code space $W$ and the compact group $\G$ generated by the
$\varphi_{i},\, i=1,\dots,N$ from \eqref{eq:G formel}, for the skew-product
shift operator $\theta:W\times\G\rightarrow W\times\G$ with $\theta(w_{1}w_{2}\dots,\,\phi):=(w_{2}w_{3}\ldots,\,\varphi_{w1}^{-1}\circ\phi)$.
Recall $\varphi_{i}$ (here the Rokhlin co-cycle) is the orthogonal
group component of $S_{i}$, $\Hl$ is the unique normalized Haar
measure on $\G$, and $\mu_{F}\left(A\right)=\H^{D}\left(F\cap A\right)/\H^{D}\left(F\right)$.
Let $P_{\Rd}$ be the projection onto the first component, and $P_{\G}$
onto the second.
\begin{fact}
$\left[W\times\G,\nu\otimes\H_{\G},\theta\right]$ is ergodic. (\cite[Proposition 5.1]{Graf95TangentialMR1363139})
\end{fact}
According to \eqref{bernoulli-meas}, $\left[W\times\G,\nu\otimes\H_{\G},\theta\right]$
induces the ergodic dynamical system $[F\times\G,\,\mu_{F}\otimes\H_{\G},\, T]$,
where the transformation $T:F\times\G\rightarrow F\times\G$ is defined
for $\mu_{F}\otimes\H_{\G}$-a.a. $\left(x,\phi\right)$ by \[
T\left(x,\phi\right)\assign\left(S_{j}^{-1}x,\varphi_{j}^{-1}\circ\phi\right)\;\mbox{if}\; x\in S_{j}(F),\; j=1,\ldots N,\]
taking into regard that $\mu_{F}(S_{i}(F)\cap S_{j}(F))=0,\, i\neq j$.
(More general references on this subject may be found, e.g., in Falconer
\cite{Falconer97TechniquesMR1449135}, Mauldin and Urbanski \cite{MauldinUrbanski03GraphDirectedMR2003772}.)
Recall we identify a.a. points in $F$ with their coding sequences.

Now we will show the curvature located on the covariant neighborhood
$A_{F}$ is covariant under the shift map. Note the transformed $\epsilon$
stays in the domain of definition: $\ep<b^{-1}d(x,(S_{x|l}J)^{c})$
implies $\ep<r_{x|l}^{-1}\ep<\ep_{0}$, since $d(x,(S_{x|i}J)^{c})=r_{x|l}\, d(P_{\Rd}T^{l}\left(x,\phi\right),J^{c})$.
From this and $A_{F}(x,\ep)\subseteq B(x,a\ep)$ we obtain for Lebesgue-a.a.
$\ep$, $\mu_{F}\otimes\H_{\G}$-a.a. $\left(x,\phi\right)\in F\times\G$,
and $l\in\mathbb{N}$ satisfying the first condition the equalities
\begin{eqnarray}
 &  & C_{k}\big(F(\ep),A_{F}(x,\ep)\times\phi(\dirs)\big)\label{eq:covariance means Birkhoff}\\
 & = & C_{k}\big(F_{x|l}(\ep),A_{F}(x,\ep)\times\phi(\dirs)\big)\nonumber \\
 & = & C_{k}\left(F_{x|l}(\ep),S_{x|l}A_{F}\left(S_{x\vert l}^{-1}\left(x\right),r_{x|l}^{-1}\ep\right)\times\phi\left(\dirs\right)\right)\nonumber \\
 & = & r_{x|l}^{k}C_{k}\left(F(r_{x|l}^{-1}\ep),A_{F}\left(S_{x\vert l}^{-1}\left(x\right),r_{x|l}^{-1}\ep\right)\times\left(\varphi_{x\vert l}^{-1}\circ\phi\left(\dirs\right)\right)\right)\nonumber \\
 & = & r_{x|l}^{k}C_{k}\left(F(r_{x|l}^{-1}\ep),A_{F}\left(P_{\Rd}T^{l}\left(x,\phi\right),r_{x|l}^{-1}\ep\right)\times P_{\G}T^{l}\left(x,\phi\right)\left(\dirs\right)\right)\,.\nonumber \end{eqnarray}
 Here we have used the locality \eqref{eq:locality-curv} of the curvature
measure $C_{k}$, the covariance of the sets $A_{F}(x,\ep,\dirs)$
(Def. \ref{def:covariant neighborhood net}), the scaling property
\eqref{eq:homogenity-curv} of $C_{k}$ under similarities, and have
expressed the inverse map of the IFS using the shift operator.

Measurability of $\left(x,\phi,\ep\right)\mapsto C_{k}\left(F(\ep),A_{F}(x,\ep)\times\phi(\dirs)\right)$
is proved the same way as in \cite[Lemma 2.3.1]{RatajZaehle10CurvatureDensitiesSelfSimilarSets}.
This works for measurability in $x,\epsilon$ for a fixed $\phi=\op{id}$.
General $\phi$ can be recovered by concatenating with the (measurable)
group action of $\phi$.

Finally, the limit of the $\Delta_{F,x,\delta}$ measures \[
\lim_{\delta\rightarrow0}\frac{1}{|\ln\delta|}\int_{\delta}^{b^{-1}d(x,J^{c})}\ep^{-k}\, C_{k}\left(F(\ep),A_{F}(x,\ep)\times\phi(\dirs)\right)\,\ep^{-1}d\ep\]
can be checked using the methods of \cite[Section 2.3]{RatajZaehle10CurvatureDensitiesSelfSimilarSets}.
The main idea is to separate the integral into chunks the Birkhoff
ergodic theorem can be applied to. The above computation shows each
chunk is a summand in the ergodic average. Their $A_{F}(x,\ep)$ has
to be replaced with our $A_{F}(x,\ep)\times\phi(\dirs)$, $d\mu_{F}\left(x\right)$
with $d\left(\mu_{F}\otimes\H_{\G}\right)\left(x,\phi\right)$, and
their $d\left(T^{i}x,J^{c}\right)$ with $d\left(P_{\Rd}T^{l}\left(x,\phi\right),J^{c}\right)$.

The proof for the variation measures $C_{k}^{\pm}$ works the same
way. (The integrability assumption of the Birkhoff theorem can be
replaced with positivity.) The extension of the limit to all $\delta\ra0$
is done by monotonicity arguments and by bounds on their $\delta\left(x,n\left(x,\delta^{\prime}\right)\right)/\delta^{\prime}$
instead of their equation (28).

\subsection{\label{sec:Global-versus-Local}Main result 2: convergence of curvature-direction
measures}

Fractal curvature-direction measures exist under slightly stronger
conditions. The local fractal curvatures $D_{C_{k}^{\Frac}|F}(\dirs)$
from Theorem \ref{curvdens} play the role of (constant) fibre measures
on $F$ with respect to the normalized Hausdorff measure $\mu_{F}$.
This is in analogy to the case of a differentiable submanifold $M$
of $\Rd$, where the (local) fibre measure on the sphere $S^{d-1}$
over a point $x\in M$ is given by a symmetric polynomial of principal
curvatures times signed unit mass atoms on the unit normals with foot-point
$x$, and where the (global) curvature-direction measure is the integral
of these fibres with respect to the intrinsic Lebesgue measure on
$M$.

Instead of adapting the proofs scattered over \cite{RatajZaehle10CurvatureDensitiesSelfSimilarSets,WinterZaehle10cmofsss,Zaehle11SelfsimRandomFractalsMR2763731,Win08MR2423952}
in a routine way, we propose a shorter proof under weaker, sharp assumptions.

Following \cite{RatajZaehle10CurvatureDensitiesSelfSimilarSets},
a Fubini argument is applied to the total curvature to move the Hausdorff
measure on $F$, which is invariant under the shift on $F$, from
inside the definition of the covariant neighborhood net, to the outermost
integral \eqref{eq:C_k <-> local version}. This connects the curvature-direction
measure with the local fibre version and the shift dynamical system
on $F$. Crucially unlike \cite{RatajZaehle10CurvatureDensitiesSelfSimilarSets},
we already work with the measure at this stage instead of its total
mass.

The group extension of the shift dynamical system accounts for the
way pieces of the attractor $F$ are rotated under the IFS. This is
not neutralized by integrating over $F$ alone, as done in the Fubini
expression, so we need a convergence result for unrotated fibre measures.
Using the Birkhoff ergodic theorem, Theorem \ref{curvdens} gave us
a result on fibres rotated in almost any way. Corollary \ref{pro:area measures converge, invariant}
improves this to unrotated ones, at the price of weak convergence. 

The dynamical system shrinks the covariant neighborhood net of each
point at its own speed, whereas the parallel set width has to be shrunk
globally. This discrepancy introduces a uniform integrability condition,
which we show to be necessary and sufficient for the convergence of
fractal curvature-direction measures to our formula (Proposition \ref{cor:sharp convergence result}).
See also Remark \ref{rem:sufficient for absolute convergence} for
an {}``easier'', sufficient condition.

Finally, Corollary \ref{cor:forDuzanPokorny} touches upon whether
the similarity dimension is always equal to the curvature scaling
dimension. If the fractal curvature vanishes trivially, then the curvature
of the parallel sets $F_{\epsilon}$ is concentrated near the boundary
$\partial J$ of the (SOSC) open set.

Recall $D_{C_{k}^{\Frac}|F}$ is defined by \eqref{curvdens-int}.
\begin{thm}
\label{thm:area measure}\label{thm:c_k converge}(Fractal curvature-direction
measure) Suppose that the self-similar set $F$ in $\Rd$ with Hausdorff
dimension $D$ satisfies the strong open set condition w.r.t. $\Int J$.
Let the system of all $A_{F}(x,\ep)$ with $x\in F$, $\ep<\ep_{0}$,
be the covariant neighborhood net given in \eqref{eq:A_F useful}
with constant $a=2c_{F}^{-1/D}$, where $c_{F}\leq1$ fulfills \eqref{D-set}
. Let $b=\max\big(2a,\ep_{0}^{-1}|J|)\big)$. 

Let $k\in\{0,1,\ldots,d\}$. If $k\le d-2$, we additionally suppose
the regularity of parallel sets, Assumption \ref{ass:Neigborhood-regularity},
and assume the function \begin{equation}
F\ni x\mapsto\frac{1}{|\ln\delta|}\int_{\delta}^{\ep_{0}}\ep^{-k}C_{k}^{\var}\big(F(\ep),B\left(x,a\epsilon\right)\times S^{d-1}\big)\,\,\ep^{-1}d\ep\label{eq:assumption fubini}\end{equation}
is uniformly $\H^{D}$-integrable for $0<\delta<\ep_{0}$.

Then the following limit exists in the weak sense: \begin{equation}
C_{k}^{\Frac}\left(F,\cdot\right)\assign\lim_{\delta\rightarrow0}\frac{1}{|\ln\delta|}\int_{\delta}^{\ep_{0}}\ep^{D-k}C_{k}\left(F(\ep),\cdot\right)\,\,\ep^{-1}d\ep.\label{eq:defining limit area measure}\end{equation}
The limit is a product and is uniquely specified by \[
C_{k}^{\Frac}(F,B\times R)=\H^{D}(F\cap B)\, D_{C_{k}^{\Frac}|F}\left(R\right)\]
for Borel sets $B\subseteq\Rd$, $\dirs\subseteq S^{d-1}$. 
\end{thm}
The assertion is true for $C_{k}^{+},C_{k}^{-},$ and $C_{k}^{\var}$,
see Proposition \ref{cor:sharp convergence result} below. Any of
the sufficient conditions in Remark \ref{rem:sufficient for absolute convergence}
implies uniform integrability \eqref{eq:assumption fubini}. The measure
localization was previously known only under these conditions.
\begin{proof}
Since the set $F\left(\epsilon_{0}\right)\times S^{d-1}$ is compact,
any continuous function on it can be approximated in norm with product
type functions. Linearity and the fact $A_{F}(x,\ep)\subseteq B\left(x,a\epsilon\right)$
then reduce the problem to Proposition \ref{cor:sharp convergence result}
below. In case $k=d$ or $k=d-1$, uniform integrability is due to
\eqref{eq:ck_var bounded} .
\end{proof}
The following corollary to Theorem \ref{curvdens} strengthens it
to hold for all elements of $\G$ instead of $\H_{\G}$-almost all,
at the price of weak convergence. Geometrically, it means $D_{C_{k}^{\Frac}|F}^{\pm}$
is the {}``distribution'' of normal rays on the fractal, weighted
by the $k$-th higher order mean curvature.
\begin{cor}
\label{pro:area measures converge, invariant}Assume the conditions
of Theorem \ref{thm:c_k converge} except uniform integrability of
\eqref{eq:assumption fubini}, and let $\pm\in\left\{ +,-\right\} $.
Write \begin{equation}
\Delta_{F,x,\delta}^{\pm}\left(\dirs\right)\assign\frac{1}{|\ln\delta|}\int_{\delta}^{d(x,J^{c})/b}\ep^{-k}C_{k}^{\pm}\big(F(\ep),A_{F}(x,\ep)\times\dirs\big)\,\,\ep^{-1}d\ep,\,\;\dirs\subseteq S^{d-1}\text{Borel}.\label{eq:local version measure}\end{equation}
Then for $\mu_{F}$-a.a. $x\in F$, the following limit exists and
does not depend on $x$: \begin{equation}
\lim_{\delta\rightarrow0}\Delta_{F,x,\delta}^{\pm}\left(g\right)=D_{C_{k}^{\Frac}|F}^{\pm}\left(g\right),\;\text{for continuous }g:S^{d-1}\ra\left[0,1\right].\label{area measure limit finalversion}\end{equation}

\end{cor}
(From here on, we denote $m\left(g\right)\assign\int g\, dm$ for
a signed measure $m$ and an integrable function $g$.) 
\begin{proof}
Firstly, $\Delta_{F,x,\delta}^{\pm}$ and $D_{C_{k}^{\Frac}|F}^{\pm}$
are Borel measures on $S^{d-1}$: Measurability of the integrand is
covered in Section \ref{sub:Extended-shift-dynamical}. To show $\sigma$-additivity
using the defining formulas \eqref{curvdens-int}, \eqref{eq:local version measure},
we combine the measure property of $C_{k}^{\pm}\left(F\left(\epsilon\right),\cdot\right)$,
and the monotone convergence lemma.

Next, we will show convergence of $\Delta_{F,x,\delta}^{\pm}\circ\phi^{-1}$
for almost all $x\in F$, $\phi\in G$. Let $\mathcal{R}$ be a countable
system of open sets that generates the topology of $S^{d-1}$. We
will assume $S^{d-1}\in\mathcal{R}$, to help with the total masses.
Set \[
\Phi\assign\left\{ \phi\in G\,:\,\lim_{\delta\ra0}\Delta_{F,x,\delta}^{\pm}\left(\phi\dirs\right)=D_{C_{k}^{\Frac}|F}^{\pm}\left(R\right)\text{ for all }\dirs\in\mathcal{R},\,\mu_{F}\text{-a.a. }x\in F\right\} .\]
$\Phi$ has full $\H_{\G}$ measure since $\mathcal{R}$ is countable,
so $\Phi$ is not empty (Theorem \ref{curvdens}). 

Consider the case $D_{C_{k}^{\Frac}|F}^{\pm}$ has finite mass. A
generator of the topology is a convergence-determining class. (\cite[Theorem 2.3]{BillingsleyConvergence99MR1700749}
also holds for sequences of finite measures if the total masses converge.)
Consequently, we have weak convergence: for all $\phi\in\Phi$, almost
all $x\in F$, and all continuous $g:S^{d-1}\ra\left[0,1\right]$,
\[
\lim_{\delta\ra0}\Delta_{F,x,\delta}^{\pm}\left(g\circ\phi^{-1}\right)=D_{C_{k}^{\Frac}|F}^{\pm}\left(g\right).\]

Since $g\circ\phi$ is also continuous, we may substitute it for $g$,
\[
\lim_{\delta\ra0}\Delta_{F,x,\delta}^{\pm}\left(g\right)=D_{C_{k}^{\Frac}|F}^{\pm}\left(g\circ\phi\right).\]

As a Haar integral, $D_{C_{k}^{\Frac}|F}^{\pm}$ is invariant under
$\phi$: \[
\lim_{\delta\ra0}\Delta_{F,x,\delta}^{\pm}\left(g\right)=D_{C_{k}^{\Frac}|F}^{\pm}\left(g\right).\]

In case of infinite $D_{C_{k}^{\Frac}|F}^{\pm}$, both it and the
map $g\mapsto\liminf_{\delta\ra0}\Delta_{F,x,\delta}^{\pm}\left(g\right)$
still are $\Phi$-invariant, and $\Phi$ is dense in $\G$. A finite
covering argument shows both are infinite on every open set $\dirs\neq\varnothing$
and (any level set of) every non-vanishing $g$, thus the assertion.
\end{proof}
This well-known fact will compensate that the supports of $\mu_{F}$
and $C_{k}\left(F\left(\epsilon\right),\cdot\right)$ are close but
not identical:
\begin{lem}
\label{lem:Max under Cesaro limit}Let $X$ be a metric space, $g:X\times\left(0,\infty\right)\ra\left[0,\infty\right)$,
$g_{0}:X\ra\left[0,\infty\right)$, and let $f$ be a continuous and
bounded function $f:X\ra\R$, and $x\in X$ such that $f\left(x\right)\neq0$.
The following are equivalent, including the existence of the limits:\end{lem}
\begin{enumerate}
\item $f\left(x\right)\lim_{\delta\ra0}\frac{1}{\left|\ln\delta\right|}\int_{\delta}^{\delta_{0}}g\left(x,\epsilon\right)\epsilon^{-1}d\epsilon=f\left(x\right)g_{0}\left(x\right)$,
\item $\lim_{\delta\ra0}\frac{1}{\left|\ln\delta\right|}\int_{\delta}^{\delta_{0}}\left(\min_{y\in B\left(x,\epsilon\right)}f\left(y\right)\right)g\left(x,\epsilon\right)\epsilon^{-1}d\epsilon=f\left(x\right)g_{0}\left(x\right)$.\end{enumerate}
\begin{proof}
Note $\min_{y\in B\left(x,\epsilon\right)}\left\vert f\left(y\right)\right\vert >0$
for sufficiently small $\epsilon$. Denote by $\omega_{f}$ the modulus
of continuity. Since $\min_{B\left(x,\epsilon\right)}f\leq f\left(x\right)\leq\min_{B\left(x,\epsilon\right)}f+\omega_{f}\left(\epsilon\right)$,
it is enough to show \[
\limsup_{\delta\ra0}\frac{1}{\left|\ln\delta\right|}\int_{\delta}^{\delta_{0}}\omega_{f}\left(\epsilon\right)g\left(x,\epsilon\right)\epsilon^{-1}d\epsilon=0.\]
 Given any $\eta>0$, there is a $\delta_{1}$ such that $\omega_{f}\left(\epsilon\right)\leq\eta\left\vert \min_{B\left(x,\epsilon\right)}f\right\vert \leq\eta\left\vert f\left(x\right)\right\vert $
whenever $\epsilon<\delta_{1}$. We split the above integral at $\delta_{1}$.
The fact $g_{0}\left(x\right)$ is finite implies $\int_{\delta_{1}}^{\delta_{0}}g\left(x,\epsilon\right)\epsilon^{-1}d\epsilon<\infty$,
so that part of the limsup becomes zero. The other part is at most
$\eta g_{0}\left(x\right)\left\vert f\left(x\right)\right\vert $. \end{proof}
\begin{prop}
\label{cor:sharp convergence result}Assume the conditions of Theorem
\ref{thm:c_k converge} except uniform integrability of \eqref{eq:assumption fubini},
and let $\pm\in\left\{ +,-\right\} $. Then for all continuous $f:\Rd\ra\left[0,\infty\right)$
and $g:S^{d-1}\ra\left[0,\infty\right)$, \[
\liminf_{\delta\rightarrow0}\frac{1}{|\ln\delta|}\int_{\delta}^{\ep_{0}}\ep^{D-k}C_{k}^{\pm}\left(F(\ep),f\otimes g\right)\,\,\ep^{-1}d\ep\geq\left(\sum_{j=1}^{N}r_{j}^{D}|\ln r_{i}|\right)^{-1}\H^{D}\left(1_{F}f\right)\, D_{C_{k}^{\Frac}|F}^{\pm}\left(g\right).\]
The following are equivalent: 
\begin{enumerate}
\item The lower limit exists as a limit, and we have equality.\label{enu:fatou strict}
\item The expression \[
\left(x,\phi\right)\mapsto\frac{1}{|\ln\delta|}\int_{\delta}^{\epsilon_{0}}\ep^{-k}C_{k}^{\pm}\big(F(\ep),1_{A_{F}(x,\ep)}\otimes\left(g\circ\phi^{-1}\right)\big)\,\,\ep^{-1}d\ep\]
is uniformly $\mu_{F}\otimes\H_{\G}$-integrable for $0<\delta<\epsilon_{0}$,
or both sides of the assertion are infinite, or $f$ is zero on $F$.\label{enu:uniformly integrable}
\end{enumerate}
\end{prop}
The integral in \eqref{enu:uniformly integrable} differs from $\Delta_{F,x,\delta}\left(g\circ\phi^{-1}\right)$
only by its upper boundary.
\begin{proof}
The left side of the assertion will be assumed finite for the entire
proof.

We will reduce the problem to Theorem \ref{curvdens} and the previous
corollary, taking into regard the relationship \begin{equation}
\ep^{D}\int_{\Rd\times S^{d-1}}f\left(z\right)g\left(n\right)\, dC_{k}^{\pm}\big(F(\ep),\left(z,n\right)\big)\geq\int_{F}\left(\min_{B\left(x,a\ep\right)}f\right)C_{k}^{\pm}\big(F(\ep),1_{A_{F}(x,\ep)}\otimes g\big)\, d\mathcal{H}^{D}\left(x\right)\label{int}\end{equation}
 for a.a. $\ep<\ep_{0}$, with equality if $f$ is a constant. The
sets $A_{F}$ were defined in \eqref{eq:A_F useful} to make this
possible: ($\dirs\subseteq S^{d-1}$ Borel)\begin{eqnarray}
 &  & \ep^{D}\int_{F(\ep)\times\dirs}f\left(z\right)\, dC_{k}^{\pm}\left(F(\ep),\left(z,n\right)\right)\label{eq:C_k <-> local version}\\
 & = & \int_{F(\ep)\times\dirs}f\left(z\right)\,\H^{D}\left(F\cap B(z,\rho_{F}(z,\ep)\right)\, dC_{k}^{\pm}\left(F(\ep),\left(z,n\right)\right)\nonumber \\
 & = & \int_{F}\int_{F(\ep)\times\id\dirs}f\left(z\right)\,{\bf 1}\left(|x-z|\le\rho_{F}(z,\ep)\right)\, dC_{k}^{\pm}\left(F(\ep),\left(z,n\right)\right)\, d\mathcal{H}^{D}(x)\nonumber \\
 & \geq & \int_{F}\left(\min_{B\left(x,a\ep\right)}f\right)\,\int_{F(\ep)\times\id\dirs}{\bf 1}\left(|x-z|\le\rho_{F}(z,\ep)\right)\, dC_{k}^{\pm}\left(F(\ep),\left(z,n\right)\right)\, d\mathcal{H}^{D}(x)\nonumber \\
 & = & \int_{F}\left(\min_{z\in B\left(x,a\ep\right)}f\left(z\right)\right)C_{k}^{\pm}\big(F(\ep),A_{F}(x,\ep)\times\dirs\big)\, d\mathcal{H}^{D}\left(x\right).\nonumber \end{eqnarray}
We were allowed to apply Fubini because the integrand is positive
and all measures $\sigma$-finite. Approximating with staircase functions
yields \eqref{int}. Then we get \begin{eqnarray}
 &  & \frac{1}{|\ln\delta|}\int_{\delta}^{\ep_{0}}\ep^{D-k}\int_{F(\ep)}f\left(z\right)g\left(n\right)\, dC_{k}^{\pm}\left(F(\ep),\left(z,n\right)\right)\,\,\ep^{-1}d\ep\nonumber \\
 & \geq & \int_{F}\frac{1}{|\ln\delta|}\int_{\delta}^{\ep_{0}}\left(\min_{B\left(x,a\ep\right)}f\right)\ep^{-k}C_{k}^{\pm}\big(F(\ep),1_{A_{F}(x,\ep)}\otimes g\big)\,\ep^{-1}d\ep\, d\H^{D}(x)\,.\label{eq:c_k <-> local, integrated eps}\end{eqnarray}

Fatou's lemma lets us take the lower limit as $\delta\rightarrow0$
under this integral. The liminf on the left side was assumed finite,
so the integrand on the right is finite almost everywhere. Next, Lemma
\ref{lem:Max under Cesaro limit} permits us to replace $\left(\min_{B\left(x,a\ep\right)}f\right)$
with $f\left(x\right)$, except in case $f\left(x\right)=0$, where
the minimum already agrees with $f\left(x\right)$. \begin{eqnarray*}
 &  & \liminf_{\delta\ra0}\frac{1}{|\ln\delta|}\int_{\delta}^{\ep_{0}}\ep^{D-k}\int_{F(\ep)\times S^{d-1}}f\left(z\right)g\left(n\right)\, dC_{k}^{\pm}\left(F(\ep),\left(z,n\right)\right)\,\,\ep^{-1}d\ep\\
 & \geq & \int_{F}\liminf_{\delta\ra0}\frac{1}{|\ln\delta|}\int_{\delta}^{\ep_{0}}\left(\min_{B\left(x,a\ep\right)}f\right)\ep^{-k}C_{k}^{\pm}\big(F(\ep),1_{A_{F}(x,\ep)}\otimes g\big)\,\ep^{-1}d\ep\, d\H^{D}(x)\\
 & = & \int_{F}\lim_{\delta\ra0}\frac{1}{|\ln\delta|}\int_{\delta}^{\ep_{0}}\left(\min_{B\left(x,a\ep\right)}f\right)\ep^{-k}C_{k}^{\pm}\big(F(\ep),1_{A_{F}(x,\ep)}\otimes g\big)\,\ep^{-1}d\ep\, d\H^{D}(x)\\
 & = & \int_{F}\lim_{\delta\ra0}\frac{1}{|\ln\delta|}f\left(x\right)\int_{\delta}^{\ep_{0}}\ep^{-k}C_{k}^{\pm}\big(F(\ep),1_{A_{F}(x,\ep)}\otimes g\big)\,\ep^{-1}d\ep\, d\H^{D}(x)\\
 & = & D_{C_{k}^{\Frac}|F}^{\pm}\left(g\right)\,\int_{F}f\left(x\right)\, d\H^{D}(x)\end{eqnarray*}
Reading backwards, Theorem~\ref{curvdens} guarantees that the lower
limit is a limit and Lemma \ref{lem:Max under Cesaro limit} indeed
applicable. We have proved the asserted inequality.

Now assume uniform integrability due to assertion \eqref{enu:uniformly integrable}.
We get equality in Fatou's lemma. In case $f$ is a constant, \eqref{eq:c_k <-> local, integrated eps}
becomes an equality, too. That proves assertion \eqref{enu:fatou strict}
first for a constant $f$. The general case follows from positivity
and equal mass. Conversely, equality in Fatou's lemma implies $L_{1}$
convergence, which in turn implies uniform integrability (\cite[Chapters 6.8 and 6.18]{Doob94MeasureTheoryMR1253752}).\end{proof}
\begin{rem}
The convergence result also holds on measurable functions $f\otimes g$
of type \[
\left(x,n\right)\mapsto f\left(x\right)1_{\dirs}\left(n\right),\:\qquad f\text{ continuous},\,\dirs=\varphi_{i}R\subseteq S^{d-1},\, i=1,\dots,N,\]
i.e. $\dirs$ is a totally invariant Borel set. (This is stronger
than convergence on all continuous functions.) Instead of setting
$\phi=\id\in\G$ in the proof by invoking Corollary \ref{pro:area measures converge, invariant},
we can exploit $g\circ\phi=g$ using Theorem \ref{curvdens} directly.
Simply replace $d\mathcal{H}^{D}(x)$ with $d\left(\mathcal{H}^{D}\otimes\H_{\G}\right)\left(x,\phi\right)$,
$\int_{F}$ with $\int_{F\times\G}$, and $\op{id}$ with $\phi$.
\end{rem}
In some degenerate cases, the similarity dimension $D$ can be an
inappropriate rescaling, and $C_{k}^{\Frac}\left(F,\cdot\right)$
can vanish trivially. One such example is the unit cube $\left[0,1\right]^{d}\subseteq\Rd$,
which is generated by $2^{d}$ similarities with contraction ratio
$1/2$ and fixed points in the corners. One of both variations can
vanish even for {}``true'' fractals, for example $C_{0}^{+,\Frac}$
of the Sierpinski gasket. (See \cite{Win08MR2423952} for further
discussion.) Such a situation can be detected by the concentration
of the curvature near $\partial J$ as $\epsilon\ra0$, i.e. away
from $\mu_{F}$-most of $F$. 
\begin{cor}
\label{cor:forDuzanPokorny}In the situation of the proposition, if
\[
\lim_{\delta\rightarrow0}\frac{1}{|\ln\delta|}\int_{\delta}^{\ep_{0}}\ep^{D-k}C_{k}^{\pm}\left(F(\ep)\right)\,\,\ep^{-1}d\epsilon=0,\]
then \[
C_{k}^{\pm}\left(F\left(\epsilon\right),\,\overline{J^{c}}\left(\epsilon c\right)^{c}\times S^{d-1}\right)=0\]
for almost all $\epsilon<\epsilon_{0}$, where \[
c\assign\max\left\{ 4\, c_{F}^{-\frac{1}{D}},\left\vert J\right\vert \H^{D}\left(F\right)^{-\frac{1}{D}}\right\} -C_{F}^{-\frac{1}{D}},\]
and $C_{F}$, $c_{F}$ are constants from \eqref{D-set}.\end{cor}
\begin{rem}
Enlarging the $\sigma$-algebra can increase the mass of the absolute
variation measure. In principle, our integrability assumption on $C_{k}\left(F\left(\ep\right),\,\cdot\right)$
as a measure in $\Rd\times S^{d-1}$ could be stronger than assumptions
on $C_{k}\left(F\left(\ep\right),\,\cdot\times S^{d-1}\right)$ in
$\Rd$. \cite{RatajZaehle10CurvatureDensitiesSelfSimilarSets} worked
with the latter. However, we know no examples of parallel sets $F\left(\ep\right)$
where this matters.\\
For example, let $K\assign\left\{ x\in\R^{2}\,:\,\left|x\right|=1\right\} $
be the unit circle. The set of outpointing normals carries the positive
$C_{0}$ mass, but they always share a foot-point with an in-pointing
normal carrying negative mass.\begin{alignat*}{5}
&\sup_{A\subseteq\R^{2}\text{Borel}}&&C_{0}\left(K,A\times S^{d-1}\right)&~=~&0,\\
&\sup_{A\subseteq\R^{2}\times S^{d-1}\text{Borel}}&&C_{0}\left(K,A\right)&~=~&1.\\
\end{alignat*} \end{rem}
\begin{fact}
(\cite[Lemma 3]{Rataj08VariationConvergenceMR2373368}) If $F\left(\ep\right)$
is a finite union of convex sets for some (and therefore all) $\epsilon>0$,
then for $k=1,\dots,d$, \[
\sup_{A\subseteq\R^{2}\text{Borel}}C_{k}\left(F\left(\ep\right),A\times S^{d-1}\right)=\sup_{A\subseteq\R^{2}\times S^{d-1}\text{Borel}}C_{k}\left(F\left(\ep\right),A\right),\]
 i.e. for almost all $x\in F\left(\ep\right)$, the normal cone over
$x$ contributes either purely positive or purely negative mass.
\end{fact}
The lemma also holds for arbitrary compact sets $F$ if $\ep$ is
large enough, see the proof of \cite[Theorem 4.1]{Zaehle11SelfsimRandomFractalsMR2763731}.

\section{Examples\label{sec:Example} }

To illustrate the directions: If $F$ is the classical Sierpinski
triangle, both $D_{C_{1}^{\Frac}|F}$ and $D_{C_{0}^{\Frac}|F}$ are
a constant multiples of their classical counterpart for an equilateral
triangle, i.e. $D_{C_{1}^{\Frac}|F}$ has equal atoms on the normals
onto the three sides, and $D_{C_{0}^{\Frac}|F}$ is a constant times
the uniform distribution on $S^{1}$. ($D_{C_{2}^{\Frac}|F}$ is always
rotation invariant for sets $F$ in $\R^{2}$.)

\begin{figure}[h]
\begin{minipage}[t][1\totalheight][c]{0.47\textwidth}%
\printornot{
\centering
\def\svgwidth{\columnwidth}
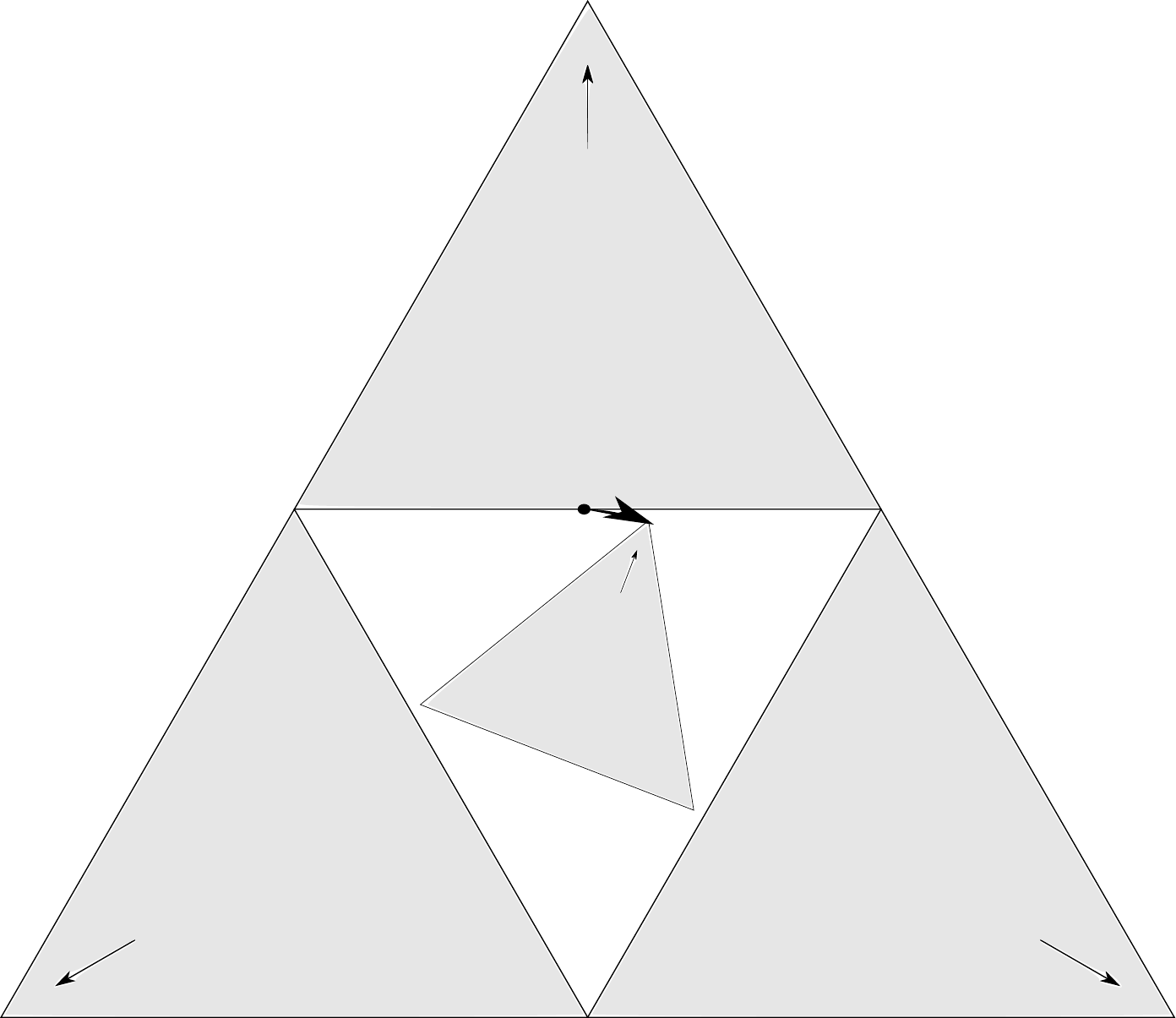
}

\caption{\selectlanguage{american}%
\label{fig:modified Sierpinski gasket IFS}IFS: images of convex hull\selectlanguage{english}
}
\end{minipage}\hfill{}%
\begin{minipage}[t][1\totalheight][c]{0.47\textwidth}%
\printornot{
\centering
\def\svgwidth{\columnwidth}
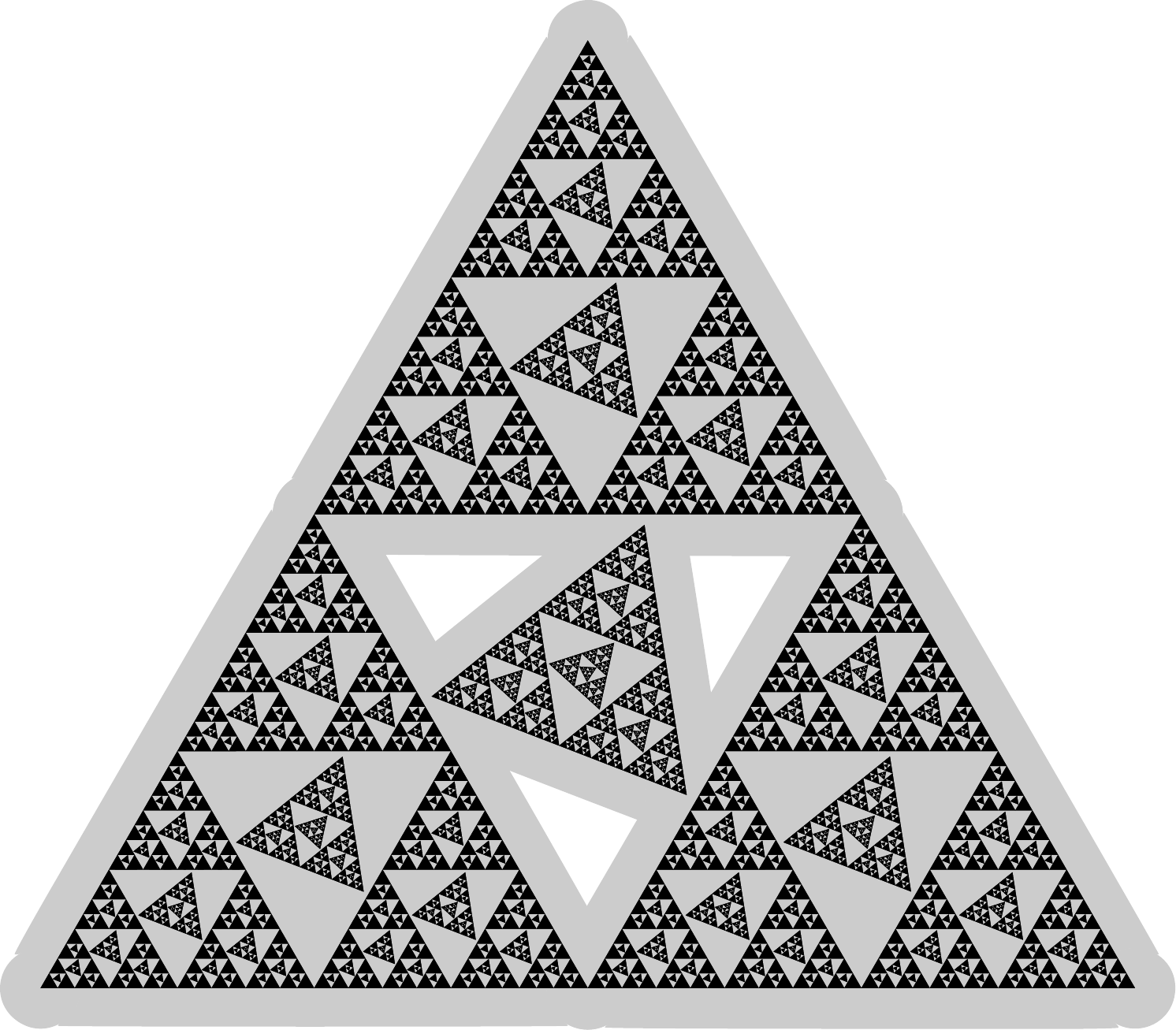
}

\caption{\selectlanguage{american}%
\label{fig:modified Sierpinski gasket parallel}Attractor $F$ and
parallel set $F\left(\epsilon\right)$ (shaded) \selectlanguage{english}
}
\end{minipage}
\end{figure}
Let $F\subseteq\R^{2}$ instead be the modified Sierpinski gasket
shown in Figure \ref{fig:modified Sierpinski gasket parallel}. It
is generated by the following IFS: $S_{i}$, $i=1,2,3$ are the same
as for the Sierpinski gasket. They contract by $r_{i}=1/2$, rotate
by $2\pi i/3\,\op{mod}2\pi$, and move the upper vertex onto the left,
right, or upper vertex, respectively. The last similarity, $S_{4}$,
contracts by $r_{4}=1/4$, rotates by $\alpha\op{mod}2\pi$, and has
its fixed point $\left(1/2,\sqrt{3}/6\right)$ in the center. ($F$
is connected only if $\alpha$ is an integer multiple of $2\pi/3$.) 

The parallel sets are regular due to polyconvexity. The Open Set Condition
is satisfied with $J\assign\op{cl}\op{conv}F$.

The compact group $\G$ generated by the rotational components is
\[
\G=\begin{cases}
S^{1} & \mbox{if }\frac{\alpha}{2\pi}\notin\mathbb{Q},\\
\left\{ k\alpha+l\frac{2\pi}{3}\op{mod}2\pi\,:\, k,l\in\Z\right\}  & \mbox{if }\frac{\alpha}{2\pi}\in\mathbb{Q}.\end{cases}\]
 \\
$D_{C_{1}^{\Frac}|F}$ exhibits only the symmetry demanded by $\G$-invariance.
Identifying $S^{1}$ and $SO\left(2\right)$ such that $0\approx\op{id}$
points downwards, $D_{C_{1}^{\Frac}|F}$ is a (positive) multiple
of the Haar measure on $\G$. (Fine) approximation sets obtained by
applying the IFS finitely often to $J$ have a $C_{1}$ measure with
the same support. By contrast, $D_{C_{0}^{\Frac}|F}$ is always a
(negative) multiple of the Haar measure on $S^{1}$.

If $\alpha=0$, both $D_{C_{0}^{\Frac}|F}$ and $D_{C_{1}^{\Frac}|F}$
agree with their classical counterpart (area measure) for a convex,
equilateral triangle.

Since $D_{C_{k}^{\Frac}|F}$ is non-null, the fractal curvature-direction
measures are the product of $D_{C_{k}^{\Frac}|F}$ with the normalized
Hausdorff measure on $F$: $C_{k}^{\Frac}\left(F,\cdot\right)=D_{C_{k}^{\Frac}|F}\otimes\mu_{F}$.

Finally, the interpretation as fractal curvatures is justified by
the fact they do not vanish trivially, i.e. $D$ is the appropriate
choice of rescaling exponent. This follows either from an explicit
computation or from \cite[Theorem 2.3.8]{Win08MR2423952}.

To illustrate the sufficient conditions for integrability (Remark
\ref{rem:sufficient for absolute convergence}): The Koch curve does
not have polyconvex parallel sets, but satisfies Winter's Strong Curvature
Bound Condition and our essential supremum \cite{Zaehle11SelfsimRandomFractalsMR2763731},
\cite[Example 5.3]{Winter10CurvBoundsArXiv10102032}. A rotated Sierpinski
carpet shows the latter is strictly weaker \cite[Example 5.2]{Winter10CurvBoundsArXiv10102032}.
The Menger sponge does not satisfy the supremum condition, but there
is an integrable dominating function that establishes uniform integrability
\cite[Example 4.4]{RatajZaehle10CurvatureDensitiesSelfSimilarSets}.
In that example, the measure version localizing $C_{k}^{\var}$ is
a new result.

\bibliographystyle{amsalpha}
\addcontentsline{toc}{section}{\refname}\bibliography{citations_bohl}

\newcommand{\etalchar}[1]{$^{#1}$}
\def\cprime{$'$} \def\lfhook#1{\setbox0=\hbox{#1}{\ooalign{\hidewidth
  \lower1.5ex\hbox{'}\hidewidth\crcr\unhbox0}}} \def\cprime{$'$}
\providecommand{\bysame}{\leavevmode\hbox to3em{\hrulefill}\thinspace}
\providecommand{\MR}{\relax\ifhmode\unskip\space\fi MR }
\providecommand{\MRhref}[2]{%
  \href{http://www.ams.org/mathscinet-getitem?mr=#1}{#2}
}
\providecommand{\href}[2]{#2}
\begin{thebibliography}{STMK{\etalchar{+}}11}

\bibitem[BG92]{BandtGraf92MR1100644}
Christoph Bandt and Siegfried Graf, \emph{Self-similar sets. {VII}. {A}
  characterization of self-similar fractals with positive {H}ausdorff measure},
  Proc. Amer. Math. Soc. \textbf{114} (1992), no.~4, 995--1001.

\bibitem[Bil99]{BillingsleyConvergence99MR1700749}
Patrick Billingsley, \emph{Convergence of probability measures}, second ed.,
  Wiley Series in Probability and Statistics: Probability and Statistics, John
  Wiley \& Sons Inc., New York, 1999, A Wiley-Interscience Publication.

\bibitem[Doo94]{Doob94MeasureTheoryMR1253752}
J.~L. Doob, \emph{Measure theory}, Graduate Texts in Mathematics, vol. 143,
  Springer-Verlag, New York, 1994.

\bibitem[Fal95]{Falconer95MinkowskiMeasurabilityMR1224615}
K.~J. Falconer, \emph{On the {M}inkowski measurability of fractals}, Proc.
  Amer. Math. Soc. \textbf{123} (1995), no.~4, 1115--1124.

\bibitem[Fal97]{Falconer97TechniquesMR1449135}
Kenneth Falconer, \emph{Techniques in fractal geometry}, John Wiley \& Sons
  Ltd., Chichester, 1997.

\bibitem[Fed59]{Federer59MR0110078}
Herbert Federer, \emph{Curvature measures}, Trans. Amer. Math. Soc. \textbf{93}
  (1959), 418--491.

\bibitem[Fer76]{Ferry76BoundaryManifoldMR0413112}
Steve Ferry, \emph{When {$\epsilon $}-boundaries are manifolds}, Fund. Math.
  \textbf{90} (1975/76), no.~3, 199--210.

\bibitem[Fu85]{Fu85TubularNeighMR816398}
Joseph Howland~Guthrie Fu, \emph{Tubular neighborhoods in {E}uclidean spaces},
  Duke Math. J. \textbf{52} (1985), no.~4, 1025--1046.

\bibitem[Gat00]{Gatzouras00LacunarityMR1694290}
Dimitris Gatzouras, \emph{Lacunarity of self-similar and stochastically
  self-similar sets}, Trans. Amer. Math. Soc. \textbf{352} (2000), no.~5,
  1953--1983.

\bibitem[Gra95]{Graf95TangentialMR1363139}
Siegfried Graf, \emph{On {B}andt's tangential distribution for self-similar
  measures}, Monatsh. Math. \textbf{120} (1995), no.~3-4, 223--246.

\bibitem[Hut81]{Hutchinson81FractalsMR625600}
John~E. Hutchinson, \emph{Fractals and self-similarity}, Indiana Univ. Math. J.
  \textbf{30} (1981), no.~5, 713--747.

\bibitem[KK10]{KombrinkKessebohmer10OneDimSelfconformal}
Marc {K}esseb{\"o}hmer and Sabrina {K}ombrink, \emph{{Fractal curvature
  measures and {M}inkowski content for one-dimensional self-conformal sets}},
  Advances of Mathematics (to appear) (2010), Arxiv 1012.5399.

\bibitem[LP93]{LapidusPomerance93MR1189091}
Michel~L. Lapidus and Carl Pomerance, \emph{The {R}iemann zeta-function and the
  one-dimensional {W}eyl-{B}erry conjecture for fractal drums}, Proc. London
  Math. Soc. (3) \textbf{66} (1993), no.~1, 41--69.

\bibitem[MU03]{MauldinUrbanski03GraphDirectedMR2003772}
R.~Daniel Mauldin and Mariusz Urba{\'n}ski, \emph{Graph directed {M}arkov
  systems}, Cambridge Tracts in Mathematics, vol. 148, Cambridge University
  Press, Cambridge, 2003, Geometry and dynamics of limit sets.

\bibitem[{Pok}11]{Pokorny11CriticalValuesSelfSimilar1101.1219}
D.~{Pokorny}, \emph{{On critical values of self-similar sets}}, Houston J. Math
  (to appear) (2011), Arxiv 1101.1219v3.

\bibitem[Rat08]{Rataj08VariationConvergenceMR2373368}
J.~Rataj, \emph{Convergence of total variation of curvature measures}, Monatsh.
  Math. \textbf{153} (2008), no.~2, 153--164.

\bibitem[RW10]{RatajWinter09MeasuresOfParallelSetsArxiv09053279}
Jan Rataj and Steffen Winter, \emph{On volume and surface area of parallel
  sets}, Indiana Univ. Math. J. \textbf{59} (2010), no.~5, 1661--1685.

\bibitem[RZ03]{RatajZaehle03NormalCyclesByApproxMR1983898}
J.~Rataj and M.~Z{\"a}hle, \emph{Normal cycles of {L}ipschitz manifolds by
  approximation with parallel sets}, Differential Geom. Appl. \textbf{19}
  (2003), no.~1, 113--126.

\bibitem[RZ05]{RatajZaehle05GeneralNormalMR2131910}
\bysame, \emph{General normal cycles and {L}ipschitz manifolds of bounded
  curvature}, Ann. Global Anal. Geom. \textbf{27} (2005), no.~2, 135--156.

\bibitem[RZ10]{RatajZaehle10CurvatureDensitiesSelfSimilarSets}
\bysame, \emph{Curvature densities of self-similar sets}, Indiana Univ. Math.
  J. (to appear) (2010), Arxiv 1009.6162.

\bibitem[Sch80]{Schneider80PolyconvexSteinerFormulaMR566443}
Rolf Schneider, \emph{Parallelmengen mit {V}ielfachheit und
  {S}teiner-{F}ormeln}, Geom. Dedicata \textbf{9} (1980), no.~1, 111--127.

\bibitem[Sch94]{Schief94SeparationMR1191872}
Andreas Schief, \emph{Separation properties for self-similar sets}, Proc. Amer.
  Math. Soc. \textbf{122} (1994), no.~1, 111--115.

\bibitem[STMK{\etalchar{+}}11]{SchroederturkMecke11MinkowskiShapeAnalysis}
G.E. Schr\"oder-Turk, W.~Mickel, S.C. Kapfer, M.A. Klatt, F.M. Schaller, M.J.F.
  Hoffmann, N.~Kleppmann, P.~Armstrong, A.~Inayat, D.~Hug, M.~Reichelsdorfer,
  W.~Peukert, W.~Schwieger, and K.~Mecke, \emph{Minkowski tensor shape analysis
  of cellular, granular and porous structures}, Advanced Materials \textbf{23}
  (2011), no.~22-23, 2535--2553.

\bibitem[Win08]{Win08MR2423952}
Steffen Winter, \emph{Curvature measures and fractals}, Dissertationes Math.
  (Rozprawy Mat.) \textbf{453} (2008), 66.

\bibitem[Win11]{Winter10CurvBoundsArXiv10102032}
\bysame, \emph{Curvature bounds for neighborhoods of self-similar sets},
  Comment. Math. Univ. Carolin. \textbf{52} (2011), no.~2, 205--226, Arxiv
  1010.2032.

\bibitem[WZ12]{WinterZaehle10cmofsss}
Steffen Winter and Martina Z{\"a}hle, \emph{Fractal curvature measures of
  self-similar sets}, Adv. in Geom. (2012), Arxiv 1007.0696.

\bibitem[Z{\"a}h86]{Zaehle86IntRepMR849863}
M.~Z{\"a}hle, \emph{Integral and current representation of {F}ederer's
  curvature measures}, Arch. Math. (Basel) \textbf{46} (1986), no.~6, 557--567.

\bibitem[Z{\"a}h11]{Zaehle11SelfsimRandomFractalsMR2763731}
\bysame, \emph{Lipschitz-{K}illing curvatures of self-similar random fractals},
  Trans. Amer. Math. Soc. \textbf{363} (2011), no.~5, 2663--2684.

\end{thebibliography}

\end{document}